\documentclass[12pt]{article}
\title{Searching for simultaneous arithmetic progressions on
elliptic curves}
\author{Irene Garc\'{\i}a--Selfa$^*$ \and Jos\'e M. Tornero\footnote{Both authors supported by FQM 218 and
BFM 2001--3207 and FEDER.}} 
\date{November, 2.004}

\begin{document}
\maketitle

\abstract{We look for elliptic curves featuring rational points
whose coordinates form two arithmetic progressions, one for each
coordinate. A constructive method for creating such curves is
shown, for lengths up to $5$.}

\vspace{.3cm}



\section{Introduction}

Let us consider an elliptic curve $E$ defined over ${\bf Q}$ by a
general Weierstrass equation
$$
Y^2 + a_1 XY + a_3 Y = X^3 + a_2 X^2 + a_4 X + a_6, \; a_i \in
{\bf Q}.
$$

\vspace{.5cm}

\noindent
{\bf Definition.--} We will say that the points $P_0,...,P_n \in E$ are in (or form an)
$x$--{\em arithmetic progression} if their $x$--coordinates are. The symmetric concept of 
$y$--arithmetic progression is defined analogously.

We will note by $S_x (E)$ and $S_y(E)$ the maximal
number of points in $x$--arithmetic progression and
$y$--arithmetic progression respectively that can be found in $E$.

\vspace{.5cm}

\noindent {\bf Remark.--} 
What we ignore about arithmetic progressions on elliptic curves is far more than
what we know. Apparently the first one in considering the problem was S.P.
Mohanty (\cite{M}) who focused on the Mordell equation $Y^2=X^3+k$
and looked for integral points forming arithmetic progressions of
difference $1$. He proved that for all these curves and for
just these progressions $S_x (E) \leq 2$ and $S_y (E) \leq 4$.

Later on, Lee and V\'elez (\cite{LV}) fixed their attention in the same family,
but they took into consideration all possible progressions. They
found infinite families of such curves verifying $S_x (E) \geq 4$
and (not simultaneously) $S_y (E) \geq 6$.

Bremner, Silverman and Tzanakis (\cite{BST}) took another quite
popular family of curves, $Y^2= X(X^2-n^2)$, and proved that for
all these curves and considering only integral points $S_x(E)
\leq 5$. They went further and proved very interesting results on
this line concerning free subgroups of rank one in arbitrary
elliptic curves. Their proofs were quite lengthy, involving very
delicate computations of local heights.

After this (although it was published earlier) Bremner (\cite{B})
carried out very clever computations in order to show that there
are infinitely many curves verifying $S_x (E) \geq 8$. Campbell
(\cite{C}) followed his lines to produce curves with eight points
in $x$--arithmetic progression and also a genus $1$ curve with
$12$ such points (unfortunately the curve was not in Weierstrass form!). Maybe
the more intriguing parts of Bremner's results are, on one side,
the numerical evidence of the fact that the length of
$x$--arithmetic progressions on elliptic curves may well not be
bounded (although Bremner himself did not risk to state such a
conjecture) and, on the other hand, the apparent connection between
long arithmetic progressions and high ranks of the Mordell-Weil
group.

\vspace{.5cm}

From these last papers it becomes clear that the main problem, when 
one deals with arithmetic progressions on elliptic curves, is that 
the number of parameters involved (if one wants to work with full 
generality) becomes unmanageable. This is why the only precise results 
known are bounded to  one--parametric families. 

\vspace{.5cm}

This paper is devoted to finding simultaneous arithmetic progressions 
on elliptic curves. The precise definition goes as follows:

\vspace{.5cm}

\noindent {\bf Definition.--} We will say that the points $P_0,...,P_n \in E$ 
are in (or form a) {\em simultaneous arithmetic progression} if their 
$x$--coordinates and its $y$--coordinates are (maybe not in the same order). 

We will note by $S_{x,y} (E)$ the maximal number of points in  simultaneous 
arithmetic progressions that can be found in $E$.

\vspace{.5cm}

\noindent {\bf Remark.--} When one looks for such progressions with more than 
three points, things start to become difficult, as the ordering of the points in both 
progressions can not coincide (see below for a precise explanation of this). 
Our initial aim was following \cite{B} and imposing the conditions with full generality
in order to narrow the search with respect to the $x$--arithmetic progression problem. 
Unfortunately it eventually became too complicated to deal with as well.
So, we took a different point of view from Bremner's (\cite{B}): we have
worked with (almost) arbitrary elliptic curves, but we have looked for a
restricted class of arithmetic progressions. With this starting point, 
we have been able to prove the following results, which were unknown
so far:

\vspace{.5cm}

\noindent {\bf Theorem 1.--} There are elliptic curves over ${\bf Q}$ with 
$S_y(E) \geq 7$.

\vspace{.5cm}

\noindent {\bf Theorem 2.--} There are elliptic curves over ${\bf Q}$ with 
$S_{x,y}(E) \geq 5$.

\section{A construction scheme}

As it is well--known (\cite{S}), any change of variables preserving
the Weierstrass form of $E$ must be of the form
$$
X' = u^2X+r, \; \; Y' = u^3Y+sX+t,
$$
so the existence (and the length) of $x$--arithmetic progressions is not
affected by changes of variables. This also implies that, up to
such a change, we can consider all the points in a certain $x$--arithmetic
progression to be integral (the $y$ case being symmetric). Therefore as an
immediate corollary of Siegel's theorem (see \cite{Siegel} for the original 
proof or \cite{S} for a modern one) $E$ cannot contain infinite $x$--arithmetic 
progressions.

Hence we are left to study if there is a universal bound,
independent of the chosen curve, for the length of arithmetic
progressions. So, assume we have a curve with an $x$--arithmetic
progression on it, say $P_0,...,P_n$, and assume $[2]P_0 \neq O$. 
Then we can take $P_0$ to $(0,0)$ and rotate the axes in order to 
take the tangent at $(0,0)$ to the line $X=0$ (warning: this may dismantle
$y$--arithmetic progressions in the original model). Then our curve must look like
$$
E: \; Y^2 + a XY + b Y = X^3 + c X^2.
$$

Now, if we make the change
$$
X \longmapsto \left( \frac{-c}{b} \right)^3 X, \; \; Y \longmapsto \left( \frac{-c}{b}
\right)^2 Y,
$$
our curve will be defined by
$$
E(a,b): \; Y^2 + a XY + b Y = X^3 - b X^2.
$$

This equation (also called Tate normal form) features two another
obvious points in $E(a,b)$ other than $(0,0)$: $(b,0)$ and
$(0,-b)$. Hence, in what follows we will make a new (strong indeed)
assumption and suppose $P_1 = (b,0)$. This implies that the
difference in our $x$--arithmetic progression must be precisely
$b$. In fact, we will actually take $P_0 = (0,-b)$ in order to
avoid the repetition of $0$ in the $y$--progression.

We must look for conditions which assure us that points $(kb,
y_k)$ appear in $E(a,b)$. Furthermore, we want these points to
form as well a $y$--arithmetic progression. Although the
subindex of the points will represent the increasing order in the
first coordinate it is obvious that, as for the second coordinate
is concerned, the subindex of a point might have nothing to do with
its position in the $y$--arithmetic progression.

In order to do that mind that if $P_k = (kb, y_k) \in E(a,b)$ then
it must hold
$$
y_k = \frac{-b(ak+1) \pm b \sqrt{(ak+1)^2+4k^2b(k-1) }}{2},
$$
so proving the existence of $P_k$ is equivalent to finding a rational
solution for the diophantine equation
$$
Z_k^2 = \left( ak +1 \right)^2 + 4k^2(k-1)b.
$$

We make now the change of variables
$$
\alpha_k = ak+1+Z_k, \; \; \beta_k = ak+1-Z_k,
$$
and so our previous equation becomes
$$
\alpha_k \beta_k + 4k^2(k-1)b = 0.
$$

Now, when we gather together the equations for $P_2,...,P_n$ we
must take into account that, for all $k$, $\alpha_k +\beta_k =
2ak+2$. The diophantine system which is equivalent to the
existence of our $x$--arithmetic progression is, therefore
$$
\begin{array}{ccccccc}
\alpha_2 \beta_2 & & & & + 16b & = & 0 \\
& \alpha_3 \beta_3 & & & + 72b & = & 0 \\
& & \ddots & & \vdots & \vdots & \vdots \\
& & & \alpha_n \beta_n & + 4n^2(n-1)b & = & 0 \\
3 \alpha_2 + 3 \beta_2 & - 2 \alpha_3 - 2 \beta_3 & & & & = & 2 \\
\vdots & & \ddots & & & \vdots & \vdots \\
n \alpha_2 + n \beta_2 & & & - 2 \alpha_n - 2 \beta_n & & = & 2(n-2)
\end{array}
$$

This can be seen as the intersection of $(n-1)$ hyperquadrics and
$(n-2)$ hyperplanes in the affine $(2n-1)$--dimensional space over
the rationals. Such a system is clearly unmanageable for, say
$n=10$ (not to say less). Observe that, as $b=0$ leads to no
progression at all, we must in fact ask for all $\alpha_i$ and
$\beta_i$ to be non--zero.

What looks particularly useful with this formulation of the
problem is that, in this context, we can write
$$
y_k = \frac{b}{2} \left( -(ak+1) \pm Z_k \right),
$$
and hence the two solutions for $y_k$ are precisely,
$-b\alpha_k/2$ and $-b\beta_k/2$. We can choose freely one of
them, as the above equations are symmetric in $\left\{ \alpha_i,
\beta_i \right\}$. This is quite useful in order to force the
existence of a simultaneous $y$--arithmetic progression.

\section{Numerical results}

Unlike the case of $x$--arithmetic progressions (for instance, in
\cite{B}), where the conditions for the existence of a
progression of length, say, $k$ are also to be filled for the
existence of a longer progression, in our case, if a set
$P_0,...,P_n$ displays a simultaneous arithmetic progression,
that does not mean, in principle, that $P_0,...,P_k$ also does 
(although we have not been able so far to find an example of this).

All the calculations in this section have been carried out with
MapleV and PARI/GP.

\vspace{.5cm}

\noindent {\bf Lenght $3$.} There are infinitely many curves with 
a simultaneous arithmetic progression of length $3$. In fact, we may even
ask both coordinates to share the same order in the progression. This 
clearly implies the three points must be collinear, and it also explains why
one cannot hope to have such examples with longer lengths.

For instance, all curves of the family
$$
E(b): Y^2 + (2b-1)XY + bY = X^3-bX^2
$$
\noindent have such a progression:
$$
\left\{  (0,-b), \, (b,0), \, (2b,b) \right\}.
$$

\vspace{.5cm}

\noindent {\bf Length $4$.} The system needed for the existence of an
$x$--progression of length $4$ is 
$$
\begin{array}{ccccc}
\alpha_2 \beta_2 & & +16b & = & 0 \\
& \alpha_3 \beta_3 & +72b & = & 0 \\
3 (\alpha_2+\beta_2) & - 2 (\alpha_3+\beta_3) & & = & 2
\end{array}
$$
with $y_2 = -\beta_2b/2$, $y_3=-\beta_3b/2$. Now, we can choose 
two values for $\beta_2$ and $\beta_3$ which will
guarantee the existence of the $y$--arithmetic progression. These values,
when substituted in the system will lead to a system of three linear
equations in $\alpha_2,\alpha_3,b$ whose matrix is
$$
\left( \begin{array}{ccc|c}
\beta_2 & 0 & 16 & 0 \\
0 & \beta_3 & 72 & 0 \\
3 & -2 & 0 & 2-3\beta_2+2\beta_3
\end{array} \right).
$$

Hence we will have a unique solution if $\beta_2 \neq \beta_3/3$, a
family of solutions if $\beta_2=-2/3,\, \beta_3=-2$ and no solutions
at all otherwise. The pair $(\beta_2 = -2/3,\beta_3=-2)$ does not
guarantee a length $4$ $y$--progression but it will appear later,
in the length $5$ study.

If we want $0$ and $-b$ to be in
the $y$--progression, there are only six possibilities for this
sequence and they are precisely
$$
\left\{ (2b,b,0,-b),(b,0,-b,-2b), (0,-b,-2b,-3b),
\left(\frac{b}{2},0,\frac{-b}{2},-b \right), \right.
$$
$$
\left. \left( 0 ,\frac{-b}{2},-b,\frac{-3b}{2} \right),
\left(0, \frac{-b}{3},\frac{-2b}{3},-b \right) \right\}.
$$

Each of them allows two possible choices for $\beta_2$ and
$\beta_3$, except the penultimate case, in which
$(\beta_2=1,\beta_3=3)$ is forbidden. The results found are shown in the
following table, except the cases $(\beta_2=-2,\beta_3=-4)$, which leads 
to a degenerate case $b=0$, and $(\beta_2=4,\beta_3=6)$, which gives 
the same curve as $(\beta_2=4,\beta_3=-2)$. 
$$
\begin{array}{cccc}
\quad \left( \beta_1, \beta_2 \right) \quad & \quad (a,b) \quad &
\; \; S_x \; \; & \; \; S_y \; \; \\ \hline \hline \\
(-4,-2) & (-5/3,-1/6) & \geq 5 & \geq 5 \\
(-2,4) & (-7/15,4/15) & \geq 5  & \geq 4 \\
(-1,1) & (-29/48,7/192) & \geq 4 & \geq 4 \\
(2/3,4/2) & (-7/9,2/27) & \geq 4 & \geq 5 \\
(1,-1) & (-5/16,1/64) & \geq 6  & \geq 7  \\
(4/3,2/3) & (-7/45,-1/270) & \geq 4 & \geq 4 \\
(3,1) & (29/96,-5/128) & \geq 4 & \geq 4 \\
(4,-2) & (1/3,1/6) & \geq 4 & \geq 5  \\
(6,4) & (25/21,-2/7) & \geq 6  & \geq 4 \\
\end{array}
$$

This search gave us the first interesting example, announced in Theorem 
1: the existence of an elliptic curve with a $y$--arithmetic progression
of length $7$, a fact not reported until now, as far as we know. We will 
look closer at this example below.

\vspace{.5cm}

\noindent {\bf Length $5$.} The system needed for the $x$--progression
is
$$
\begin{array}{cccccc}
\alpha_2 \beta_2 & & & +16b & = & 0 \\
& \alpha_3 \beta_3 & & +72b & = & 0 \\
& & \alpha_4 \beta_4 & +192b & = & 0 \\
3 (\alpha_2+\beta_2) & - 2 (\alpha_3+\beta_3) & & & = & 2 \\
4 (\alpha_2 + \beta_2) & & -2 (\alpha_4+\beta_4) & & = & 4
\end{array}
$$
where $y_j = -\beta_jb/2$, for $j=2,3,4$. Observe now that a blind
choice of $(\beta_2,\beta_3,\beta_4)$ as in the previous case will
lead, in general, to an incompatible system of equations.

In fact, the rank of the coefficients matrix is $4$, except in the case
$\beta_2= \beta_3/3=\beta_4/6$ but, in this case the system is
incompatible. However, these relations will prove useful later on.

There are ten possible progressions of length $5$ containing
both $0$ and $-b$, each of them permitting six different choices
for the triple $(\beta_2,\beta_3,\beta_4)$, which are the
permutations of a single choice. From this $60$ cases only two led
to a compatible system. This two cases, shown below, prove therefore Theorem 2.


\vspace{.3cm}

\noindent \fbox{$\beta_2=-4, \, \beta_3=-2, \, \beta_4=-6$}


\vspace{.3cm}

This triple gives the curve $E(-5/3,-1/6)$ with the following
points lying on it:
$$
\left\{ \left( 0,\frac{1}{6} \right), \left( \frac{-1}{6}, 0 \right),
\left( \frac{-2}{6}, \frac{-2}{6} \right), \left( \frac{-3}{6},
\frac{-1}{6} \right), \left( \frac{-4}{6}, \frac{-3}{6} \right) \right\},
$$
which form a simultaneous arithmetic progression of length $5$.

Note that, in this case, $\alpha_2=-2/3$, hence it can be
considered as a subcase of the infinite family of $x$--arithmetic
progressions found with lenght $4$.


\vspace{.3cm}

\noindent \fbox{$\beta_2=1, \, \beta_3=-1, \, \beta_4=-2$}


\vspace{.3cm}

The resulting curve is $E(-5/16,1/64)$ with the following points:
$$
\left\{ \left( 0,\frac{-2}{128} \right), \left( \frac{1}{64}, 0 \right),
\left( \frac{2}{64}, \frac{-1}{128} \right), \left( \frac{3}{64},
\frac{1}{128} \right), \left( \frac{4}{64}, \frac{2}{128} \right) \right\},
$$
which form a simultaneous arithmetic progression (note that this curve also
appeared in the previous case). Furthermore, other points lying on the curve are
$$
\left( \frac{1}{8}, \frac{-4}{128} \right)
\left( \frac{-1}{32}, \frac{-3}{128} \right), \;
\left( \frac{5}{64}, \frac{-1}{64} \right),
$$
hence as noted above $S_x \geq 6$, $S_y \geq 7$, although there are no 
simultaneous progressions in the curve of length $6$. 

\vspace{.15cm}

\noindent {\bf Remark.--} Please note that, as we mentioned at the beginning of 
the section, both examples of simultaneous arithmetic progressions of length $5$ also have 
arithmetic progressions of length $4$. This is also the case with the nine 
examples of length $4$ found; however, these data seem to us not enough for
conjecturing that this holds in general.

\vspace{.15cm}

Back to our search, another reasonable way of constructing simultaneous arithmetic 
progressions of length $5$ seems to be using $0$ and $-b(a+1)$ as
terms of our progression and again we have $60$ possibilities, all of them useless.
Most cases can be discarded using one of the following arguments, explained with two examples.


\vspace{.3cm}

\noindent \fbox{$\beta_2=4(a+1),\beta_3=6(a+1),\beta_4=8(a+1)$}


\vspace{.3cm}

If this choice made sense, we would have the $y$--progression
$$
\{ 0,-b(a+1),-2b(a+1),-3b(a+1),-4b(a+1) \}.
$$

As we said above, there are $5$ more possibilities in order to
assure this sequence, permuting the values of $\beta_2$, $\beta_3$
and $\beta_4$, but this one suits us well as illustration. Our system
matrix is now
$$
\left( \begin{array}{ccccc}
-4(a+1) & 0 & 0 & 16 & 0 \\
0 & -6(a+1) & 0 & 72 & 0 \\
0 & 0 & -8(a+1) & 192 & 0 \\
3 & -2 & 0 & 0 & 2 \\
4 & 0 & -2 & 0 & 4 \\
\end{array} \right)
$$
with determinant $-3072(a+1)^2$. Hence we must take $a=-1$ if we want the
system to have solution, but this case is degenerate.


\vspace{.3cm}

\noindent \fbox{$\beta_2=6(a+1),\beta_3=4(a+1),\beta_4=8(a+1)$}


\vspace{.3cm}

This case is also impossible, but for different reasons: it has a
matrix with determinant $7168(a+1)^2(4a+5)$. The choice
$a=-5/4$ leads to a compatible system. But remember that it also must
hold
$$
a = \frac{\alpha_i + \beta_i - 2}{2i}, \; i=2,3,4,
$$
which is not true in this case. In fact, this extra condition
annihilates the advantage of making a parameter choice, because
the value of the parameter cannot be truly arbitrary.

\section{Final remarks}

The above arguments, when applied to a length $n$ simultaneous
arithmetic progression lead to $n(n-1)/2$ progressions with
$(n-2)!$ possible value choices for each one. That is, $n!/2$
systems have to be checked. We have done the $360$ calculations
for $n=6$, obtaining no simultaneous progressions. Of course, that
does not mean that the case $n=7$ will also be unsuccessful, as we
noted before, but most probably this kind of search will not produce 
further results.

As a final comment, we would like to stress possible ways for
expanding the results in this paper:

\begin{enumerate}
\item[(a)] The first obvious thing to do it is enlarging the scope of
the progressions considered. Full generality, although desirable, may
be too messy for dealing with, at least with the current state of the art. A simpler
way could be, for instance, considering $x$--progressions in which
$b$ appears, although not necessarily as the first term. However, this has also 
proved to be unsuccessful for $n=6$ (at a cost of around 72 CPU hours).

\item[(b)] While looking for $y$--progressions seems difficult, the
formulation used here for the existence of $x$--progressions of given
length can be found useful. In a future paper, we expect to explore this
presentation of the problem, with the help of elimination theory and Gr\"obner bases.
\end{enumerate}

\vspace{.3cm}

Irene Garc\'{\i}a--Selfa (E--mail: \verb|igselfa@us.es|).

Dep. de \'Algebra, Fac. de Matem\'aticas, Universidad de Sevilla.

Apdo. 1160, 41080 Sevilla (Spain).

\vspace{.3cm}

Jos\'e M. Tornero (E--mail: \verb|tornero@algebra.us.es|).

Dep. de \'Algebra, Fac. de Matem\'aticas, Universidad de Sevilla.

Apdo. 1160, 41080 Sevilla (Spain).

\end{document}